\documentclass[runningheads,a4paper]{llncs}

\usepackage{amssymb}
\usepackage{graphicx}

\usepackage{url}
\newcommand{\keywords}[1]{\par\addvspace\baselineskip
\noindent\keywordname\enspace\ignorespaces#1}

\newcommand{\1}{1 \hspace*{-0.2ex}\rule{0.10ex}{1.5ex}\hspace{0.2ex}}

\begin{document}

\mainmatter  % start of an individual contribution

\title{A review on global sensitivity analysis methods}

\titlerunning{Global sensitivity analysis}

\author{Bertrand Iooss\inst{1}\inst{2} \and Paul Lema\^{\i}tre\inst{1}\inst{3}}
\authorrunning{Global sensitivity analysis}

\institute{EDF R\&D, 6 Quai Watier - 78401 Chatou, France
\and Institut de Math\'ematiques de Toulouse, 118 route de Narbonne - 31062 Toulouse, France
\and INRIA Sud-Ouest, 351 cours de la lib\'eration - 33405 Talence, France}

\toctitle{Global sensitivity analysis}
\tocauthor{Bertrand Iooss and Paul Lema\^{\i}tre}
\maketitle

%%%%%%%%%%%%%%%%%%%
\begin{abstract}

This chapter makes a review, in a complete methodological framework, of various global sensitivity analysis methods of model output.
Numerous statistical and probabilistic tools (regression, smoothing, tests, statistical learning, Monte Carlo, \ldots) aim at determining the model input variables which mostly contribute to an interest quantity depending on model output.
This quantity can be for instance the variance of an output variable.
Three kinds of methods are distinguished: the screening (coarse sorting of the most influential inputs among a large number), the measures of importance (quantitative sensitivity indices) and the deep exploration of the model behaviour (measuring the effects of inputs on their all variation range).
A progressive application methodology is illustrated on a scholar application.
A synthesis is given to place every method according to several axes, mainly the cost in number of model evaluations, the model complexity and the nature of brought information.

\end{abstract}

%%%%%%%%%%%%%%%%%%%%%%%
\section{Introduction}

While building and using numerical simulation models, 
Sensitivity Analysis (SA) methods are invaluable tools. They allow to study
how the uncertainty in the output of a model can
be apportioned to different sources of uncertainty in the model input (Saltelli et al. \cite{salcha00}).
It may be used to determine the most contributing input variables
to an output behavior as the non-influential
inputs, or ascertain some interaction effects within the model. The
objectives of SA are numerous; one can mention model verification and understanding,
model simplifying and factor prioritization.  Finally, the SA is an aid in the validation of a computer code, guidance 
research efforts, or the justification in terms of system design safety. 

There are many application examples, for instance Makowski et
al. \cite{maknau06} analyze, for a crop model prediction,
the contribution of $13$ genetic parameters on the variance of two outputs.
Another example is given in the work of Lefebvre et al. \cite{lefrob10}
where the aim of SA is to determine the most influential input among
a large number (around $30$), for an aircraft infrared signature simulation
model.
In nuclear engineering field, Auder et al. \cite{auddec12} study the influential inputs on thermohydraulical phenomena occuring during an accidental scenario, while Iooss et al. \cite{ioovan06} and Volkova et al. \cite{volioo08} consider the environmental assessment of industrial facilities.
 
The first historical approach to SA is known as
the local approach. The impact of small input perturbations 
on the model ouput is studied. These small perturbations occur around nominal
values (the mean of a random variable for instance). 
This deterministic approach consists in calculating or estimating the partial derivatives of the model at a specific point. 
The use of adjoint-based methods
allows to process models with a large number of input variables. Such approaches
are commonly used in solving large environmental systems as in climate modeling,
oceanography, hydrology, etc. (Cacuci \cite{cac81a}, Castaings et al. \cite{casdar09}).

From the late 1980s, to overcome the limitations of local methods (linearity and normality assumptions, local variations), a new class of methods has been developed
in a statistical framework. In contrast to local sensivity analysis, it is referred to as ``global sensitivity analysis''
because it considers the whole  variation range
of the inputs (Saltelli et al. \cite{salcha00}). 
Numerical model users and modelers have shown large interests in these tools which take full advantages of the advent on computing materials and numerical methods (see Helton \cite{hel93}, de Rocquigny et al. \cite{derdev08} and Faivre et al. \cite{faiioo13} for 
industrial and environmental applications). Saltelli et al. \cite{saltar04} and Pappenberger et al. \cite{paprat10} emphasized the need to specify clearly
the objectives of a study before making a SA. 
These objectives may include:
\begin{itemize}
\item identify and prioritize the most influential inputs,
\item identify non-influential inputs in order to fix them to nominal values,
\item map the output behavior in function of the inputs by focusing on
a specific domain of inputs if necessary,
\item calibrate some model inputs using some available information (real output observations, constraints, etc.).
\end{itemize}

With respect to such objectives, first syntheses on the subject of SA were developed
(Kleijnen \cite{kle97}, Frey and Patil \cite{frepat02}, Helton et al. \cite{heljoh06}, Badea and Bolado \cite{badbol08}, de Rocquigny et al.
\cite{derdev08}, Pappenberger et al. \cite{paprat10}). Unfortunately, between heuristics, graphical tools,
design of experiments theory, Monte Carlo techniques, statistical learning methods, etc., beginners and non-specialist users can be found quickly lost
on the choice of the most suitable methods for their problem. The aim of this chapter is to provide
an educational synthesis of SA methods inside an applicative methodological framework. 

The model input vector is denoted $\mathbf{X}=(X_1,\ldots,X_d) \in \mathbb{R}^d$. 
For the sake of simplicity, we restrict the study to a
scalar output $Y \in \mathbb{R}$ of the computer code (also called ``model'')  $f(\cdot)$:
\begin{equation}
Y = f(\mathbf{X})\;.
\end{equation}
In the probabilistic setting, $\mathbf{X}$ is a random vector defined
by a probability distribution and $Y$ is a random variable. 
In the following, the inputs  $X_i$ ($i=1\ldots d$) are assumed to be independent.
More advanced works, listed in the last section, take into account the dependence between
components of $\mathbf{X}$ (see Kurowicka and Cooke \cite{kurcoo06} for an
introduction to this issue).
Finally, this review focuses on the SA with respect to the global variability of the model output, usually measured
 by its variance.

\begin{figure}[!ht]
\begin{center}
    \includegraphics[width=8cm,height=6cm]{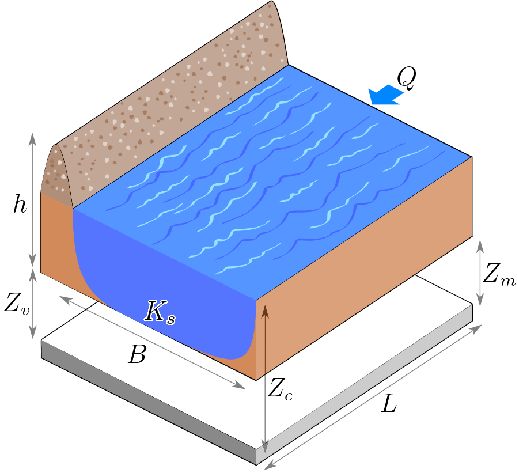} 
\end{center}
\caption{The flood example: simplified model of a river.}\label{fig:crues}
\end{figure}

All along this chapter, we illustrate our discussion with a simple application model that simulates the height of a river and compares it to the height of a dyke that protects industrial facilities (Figure \ref{fig:crues}). When the river height exceeds the one of the dyke, flooding occurs. This academic model is used as a pedagogical example in de Rocquigny \cite{der06a} and Iooss \cite{ioo11}.
The model is based on a crude simplification of the 1D hydro-dynamical equations of SaintVenant under the assumptions of uniform and constant flowrate and large rectangular sections. It consists of an equation that involves the characteristics of the river stretch:
\begin{equation}\label{eq:cruesS}
S = Z_v + H -H_d - C_b \quad \mbox{with} \quad H = \left(\frac{Q}{BK_s \sqrt{\frac{Z_m-Z_v}{L} }} \right)^{0.6},
\end{equation}
where $S$ is the maximal annual overflow (in meters), $H$ is the maximal annual height of the river (in meters) and the other variables ($8$ inputs) are defined in Table \ref{tab:factors} with their  probability distribution.
Among the input variables of the model, $H_d$ is a design parameter.
Its variation range corresponds to a design domain.
The randomness of the other variables is due to their spatio-temporal variability, our ignorance of their true value or some inaccuracies of their estimation. We suppose that the input variables are independent.
\begin{table}[!ht]
  \begin{center}
   \begin{tabular}{lccc}
Input & Description & Unit & Probability distribution \\
   \hline
 $Q$ & Maximal annual flowrate & m$^3$/s & Truncated Gumbel ${\mathcal G}(1013, 558)$ on $[500 , 3000 ]$ \\
 $K_s$ & Strickler coefficient & - & Truncated normal ${\mathcal N}(30, 8)$ on $[15 , +\infty [$ \\
 $Z_v$ & River downstream level & m & Triangular  ${\mathcal T}(49, 50, 51)$ \\
 $Z_m$ & River upstream level  & m  & Triangular  ${\mathcal T}(54, 55, 56)$  \\
 $H_d$ & Dyke height & m &  Uniform ${\mathcal U}[7, 9]$ \\
 $C_b$ & Bank level  & m & Triangular  ${\mathcal T}(55, 55.5, 56)$ \\
 $L$ & Length of the river stretch  & m &  Triangular  ${\mathcal T}(4990, 5000, 5010)$ \\
 $B$ & River width  & m &  Triangular  ${\mathcal T}(295, 300, 305)$ \\
\hline
    \end{tabular}
    \caption{Input variables of the flood model and their probability distributions.}
    \label{tab:factors}
  \end{center}
\end{table}

 We also consider another model output: the associated cost (in million euros) of the dyke,
\begin{equation}\label{eq:cruesCp}
C_p = \1_{S>0} +  \left[0.2 + 0.8\left( 1-\exp^{-\frac{1000}{S^4}}\right) \right]\1_{S \leq 0} + \frac{1}{20}\left(H_d \1_{H_d>8} + 8 \1_{H_d \leq 8} \right),
\end{equation}
with $\1_{A}(x)$ the indicator function which is equal to 1 for $x \in A$ and 0 otherwise.
In this equation, the first term represents the cost due to a flooding ($S>0$) which is $1$ million euros, the second term corresponds to the cost of the dyke maintenance ($S \leq 0$) and the third term is the investment cost related to the construction of the dyke. The latter cost is constant for a height of dyke less than $8$ m and is growing proportionally with respect to the dyke height otherwise.

The following section discusses the so-called screening methods, which are qualitative methods for studying sensitivities on
 models containing several tens of input variables. The most used quantitative measures of influence
are described in the third section.
The fourth section deals with more  advanced tools, which aim to provide a subtle exploration of the model output behavior. 
Finally, a conclusion provides a classification of these
methods and a flowchart for practitioners.
It also discusses some open problems in SA.

%%%%%%%%%%%%%%%%%%%%%%%
\section{Screening techniques}

Screening methods are based on a discretization of the inputs in levels,
allowing a fast exploration of the code behaviour. These methods
are adapted to a large number of inputs; practice has often shown that
only a small number of inputs are influential.  The aim
of this type of method is to identify the non-influential inputs with
a small number of model calls while making realistic hypotheses on the model complexity. The model is therefore simplified before
using other SA methods, more subtle but more costly.

The most engineering-used screening method is based on the so-called ``One At a Time'' (OAT) design, where each input is varied while fixing the others (see Saltelli and Annoni \cite{salann10} for a critique of this basic method).
In this section, the choice has been made to present the Morris method \cite{mor91}, which is the most complete and most costly one.
However, when the number of experiments has to be smaller than the number of inputs, one can quote the usefulness of the supersaturated design (Lin \cite{lin93}), the screening by groups (Dean and Lewis \cite{dealew06}) and the sequential bifurcation method (Bettonvil and Kleijnen \cite{betkle96}).
When the number of experiments is of the same order than the number of inputs, the classical theory of experimental design applies (Montgomery \cite{mon04}) for example with the so-called factorial fractional design.

The method of Morris allows to classify the inputs in three groups:
inputs having negligible effects, inputs having large linear effects without interactions
and inputs having large non-linear and/or interaction effects. The method
consists in discretizing the input space for each variable, then performing
a given number of OAT design. Such designs of experiments are randomly
choosen in the input space, and the variation direction is also random.
The repetition of these steps allows the estimation of elementary
effects for each input. From these effects, sensitivity indices are derived. 

Let us denote $r$ the number of OAT designs (Saltelli et al.
\cite{saltar04} propose to set parameter $r$ between
$4$ and $10$). Let us discretize the input space in a $d-$dimensionnal
grid with $n$ levels by input. Let us denote $E_{j}^{(i)}$ the elementary
effect of the $j-$th variable obtained at the $i-$th repetition,
defined as:
\begin{equation}
E_{j}^{(i)}=\frac{f(\mathbf{X}^{(i)}+\triangle e_{j})-f(\mathbf{X}^{(i)})}{\triangle}
\end{equation}
where $\triangle$ is a predetermined multiple of $\frac{1}{(n-1)}$
and $e_{j}$ a vector of the canonical base. Indices are obtained
as follows:
\begin{itemize}
\item $\mu_{j}^{*}={\displaystyle \frac{1}{r}\sum_{i=1}^{r}|}E_{j}^{(i)}|$
(mean of the absolute value of the elementary effects),
\item $\sigma_{j}=\sqrt{\displaystyle \frac{1}{r}\sum_{i=1}^{r} \left( E_{j}^{(i)} - \frac{1}{r} \sum_{i=1}^{r} E_{j}^{(i)} \right)^2}$
(standard deviation of the elementary effects).
\end{itemize}
The interpretation of the indices is the following:
\begin{itemize}
\item $\mu_{j}^{*}$ is a measure of influence of the $j-$th input on the
output. The larger $\mu_{j}^{*}$ is, the more the $j-$th input contributes
to the dispersion of the output.
\item $\sigma_{j}$ is a measure of non-linear and/or interaction effects of the $j-$th input.
If $\sigma_{j}$ is small, elementary effects have low variations on the support
of the input. Thus the effect of a perturbation is the same all along
the support, suggesting a linear relationship between the studied
input and the output. On the other hand, the larger $\sigma_{j}$
is, the less likely the linearity hypothesis is. Thus a variable with
a large $\sigma_{j}$ will be considered having non-linear effects,
or being implied in an interaction with at least one other variable.
\end{itemize}
Then, a graph linking $\mu_{j}^{*}$ and $\sigma_{j}$ allows to distinguish
the 3 groups. 

\begin{figure}[!ht]
        \centering
    \includegraphics[width=5.5cm,height=5.5cm]{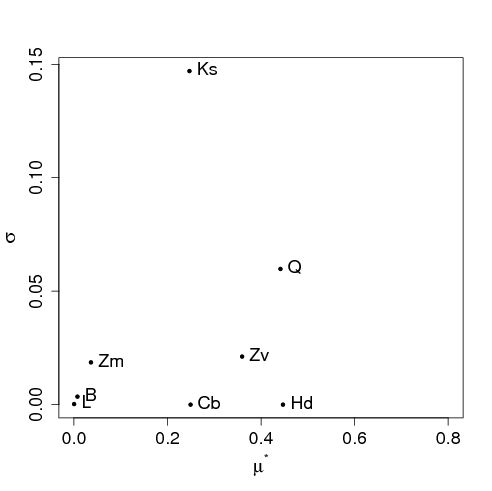} \includegraphics[width=5.5cm,height=5.5cm]{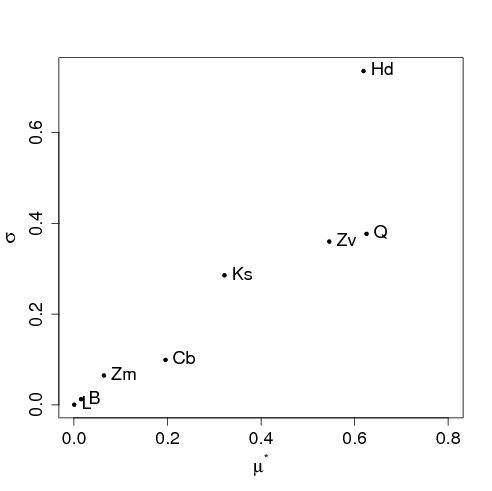} 
\caption{Results of Morris method ($r=5$ with $4$ levels): outputs $S$ (left) and $C_p$ (right).}\label{fig:crues_mor}
\end{figure}

Morris method is applied on the flood example (Eqs. (\ref{eq:cruesS}) and (\ref{eq:cruesCp})) with $r=5$ repetitions, which require $n=r(p+1)=45$ model calls.
Figure \ref{fig:crues_mor} plots results on the graph $(\mu_j^*,\sigma_j)$.
This vizualisation allows to make the following discussion:
\begin{itemize}
\item output  $S$: inputs $K_s$, $Z_v$, $Q$, $C_b$ et $H_d$ are influent, while other inputs have no effects.
In addition, the model output linearly depends on the inputs and there is no input interaction (because $\sigma_j \ll \mu_j^* \forall j$).
\item output $C_p$: inputs $H_d$, $Q$, $Z_v$ et $K_s$ have strong influence with non-linear and/or interaction effects (because $\sigma_j$ and $\mu_j^*$ have the same order of magnitude).
$C_b$ has an average influence while the other inputs have no influence.
\end{itemize}
Finally, after this screening phase, we have identified that three inputs ($L$, $B$ and $Z_m$) have no influence on the two model outputs
In the following, we fix these three inputs to their nominal values (which are the modes of their respective triangular distributions).

%%%%%%%%%%%%%%%%%%%%%%%
\section{Importance measures}

%%%%%%%%%%%%%
\subsection{Methods based on the analysis of linear models}\label{sec:regr}

If a sample of inputs and outputs $(\mathbf{X}^n,\mathbf{Y}^n)=\left(X_1^{(i)},\ldots,X_d^{(i)},Y_i\right)_{i=1..n}$ is available, it is possible to fit a linear
model explaining the behaviour of $Y$ given the values of $\mathbf{X}$, provided that the sample size $n$ is sufficiently large (at least $n > d$). 
Some global sensitivity measures defined through
the study of the fitted model are presented in the
following. 
Main indices are:
\begin{itemize}
\item {Pearson correlation coefficient:}
\begin{equation}
\rho(X_{j},Y)=\frac{\sum_{i=1}^{N}(X_{j}^{(i)}-\mathbb{E}(X_{j}))(Y_i-\mathbb{E}(Y))}
{\sqrt{\displaystyle \sum_{i=1}^{N}\left(X_{j}^{(i)}-\mathbb{E}(X_{j})\right)^{2}}\sqrt{\displaystyle \sum_{i=1}^{N}\left(Y_i-\mathbb{E}(Y)\right)^{2}}}\;.
\end{equation}
It can be seen as a linearity measure between variable $X_{j}$ and
output $Y$. It equals $1$ or $-1$ if the tested input variable
has a linear relationship with the output. If $X_{j}$ and $Y$ are
independents, the index equals $0$.
\\
\item {Standard Regression Coefficient (SRC):}
\begin{equation}
\mbox{SRC}_{j}=\beta_{j}\sqrt{\frac{\mbox{Var}(X_{j})}{\mbox{Var}(Y)}}
\end{equation}
where $\beta_{j}$ is the linear regression coefficient associated
to $X_{j}$. $\mbox{SRC}_{j}^2$ represents a share of variance if the linearity hypothesis
is confirmed.
\\
\item {Partial Correlation Coefficient (PCC):}
\begin{equation}
\mbox{PCC}_{j}=\rho(X_{j}-\widehat{X_{-j}},Y-\widehat{Y_{-j}})
\end{equation}
where $\widehat{X_{-j}}$ is the prediction of the linear model, expressing
$X_{j}$ with respect to the other inputs and $\widehat{Y_{-j}}$ is
the prediction of the linear model where $X_{j}$ is absent. PCC measures
the sensitivity of $Y$ to $X_{j}$ when the effects of the other
inputs have been canceled. 
\end{itemize}
The estimation of these sensitivity indices is subject to an uncertainty estimation, due to the limited size of the sample.
Analytical formulas can be applied in order to estimate this uncertainty (Christensen \cite{chr90}).

These three indices are based on a linear relationship between the output and the inputs.
Statistical techniques allow to confirm the linear hypothesis, as the classical coefficient of determination $R^2$ and the predictivity coefficient $Q^2$ (also called the Nash-Sutcliffe model efficiency):
\begin{equation}\label{eqQ2}
Q^2 = 1 - \frac{\sum_{i=1}^{m} [ Y^p_i-\widehat{Y}(\mathbf{X}^{p(i)}) ]^2}{\sum_{i=1}^{m} [Y^p_i-\mathbb{E}(Y^p)]^2}
\end{equation}
where $(\mathbf{X}^{p(i)},Y^p_i)_{i=1..m}$ is a $m$-size test sample of inputs-output (not used for the model fitting) and $\widehat{Y}(\cdot)$ is the predictor of the linear regression model.
The value of $Q^2$ corresponds to the percentage of output variability explained by the linear regression model (a value equals to $1$ means a perfect fit).
If the input variables are independent, each SRC$_j^2$ expresses the part of output variance explained by the input $X_j$.

If the linear hypothesis is contradicted, one can use the same three importance measures (correlation coefficient, SRC and PCC) than previously using a rank transformation (Saltelli et al. \cite{salcha00}).
The sample $(\mathbf{X}^n,\mathbf{Y}^n)$ is transformed into a sample $(\mathbf{R}_{X}^n,\mathbf{R}_{Y}^n)$ by replacing the values by their ranks in each column of the matrix. 
As importance measures, it gives the Spearman correlation coefficient $\rho^S$, the Standardized Rank Regression Coefficient (SRRC) and the Partial Rank Correlation Coefficient (PRCC).
Of course, monotony hypothesis has to be validated as in the previous case, with the determination coefficient of the ranks ($R^{2*}$) and the predictivity coefficient of the ranks ($Q^{2*}$).
%, par exemple à l'aide du coefficient de détermination $R^{2*}$ et du coefficient de prédictivité $Q_2^*$ associé à la régression linéaire sur les rangs.  

These linear and rank-based measures are part of the so-called sampling-based global sensitivity analysis method.
This has been deeply studied by Helton and Davis \cite{heldav03} who have shown the interest to use a Latin Hypercube Sample (Mc Kay et al. \cite{mckbec79}) in place of a Monte Carlo sample, in order to increase the accuracy of the sensitivity indices.

These methods are now applied on the flood example (Eqs. (\ref{eq:cruesS}) and (\ref{eq:cruesCp})) with the $d=5$ inputs that have been identified as influent in the previous screening exercise.
A Monte Carlo sample of size $n=100$ gives $100$ model evaluations.
Results are the following:
\begin{itemize}
	\item output $S$:\\
	SRC$^2(Q)=0.28$; SRC$^2(K_s)=0.12$; SRC$^2(Z_v)=0.15$; SRC$^2(H_d)=0.26$; SRC$^2(C_b)=0.03$ with $R^2=0.98$;\\
	SRRC$^2(Q)=0.27$; SRRC$^2(K_s)=0.12$; SRRC$^2(Z_v)=0.13$; SRRC$^2(H_d)=0.26$; SRRC$^2(C_b)=0.02$ with $R^{2*}=0.95$;
\item output $C_p$:\\
	SRC$^2(Q)=0.25$; SRC$^2(K_s)=0.16$; SRC$^2(Z_v)=0.18$; SRC$^2(H_d)=0.00$; SRC$^2(C_b)=0.07$  with $R^2=0.70$;\\
	SRRC$^2(Q)=0.26$; SRRC$^2(K_s)=0.19$; SRRC$^2(Z_v)=0.18$; SRRC$^2(H_d)=0.06$; SRRC$^2(C_b)=0.03$  with $R^{2*}=0.73$.
\end{itemize}
For the output $S$, $R^2$ is close to one, which shows a good fit of linear model on the data.
Analysis of regression residuals confirms this result.
Variance-based sensitivity indices are given using SRC$^2$.
For the output $C_p$, $R^2$ and $R^{2*}$ are not close to one, showing that the relation is neither linear nor monotonic.
SRC$^2$ and SRRC$^2$ indices can be used in a coarse approximation, knowing that it remains  $30\%$ of non-explained variance.
However, using another Monte Carlo sample, sensitivity indices values can be noticeably different.
Increasing the precision of these sensitivity indices would require a large increase of the sample size.

%%%%%%%%%%%%%%%%%%%%%%
\subsection{Functional decomposition of variance: Sobol' indices}\label{sec:sobol}

When the model is non-linear and non-monotonic,  the decomposition
of the output variance is still defined and can be used for SA. Let us have $f(.)$ a square-integrable function,
defined on the unit hypercube $[0,1]^{d}$. It is possible to
represent this function as a sum of elementary functions (Hoeffding \cite{hoe48}): 
\begin{equation}
f(\mathbf{X})=f_{0}+\sum_{i=1}^{d}f_{i}(X_{i})+\sum_{i<j}^{d}f_{ij}(X_{i},X_{j})+\cdots+f_{12...d}(\mathbf{X})\;.
\end{equation}
This expansion is unique under conditions (Sobol \cite{sob93}):
\[
\int_{0}^{1}{\displaystyle f_{i_{1}...i_{s}}(x_{i_{1}},...,x_{i_{s}})dx_{i_{k}}=0}\;, 1\leq k\leq s,\; \left\{ i_{1},...,i_{s}\right\} \subseteq\left\{ 1,...,d\right\} .
\]
This implies that $f_{0}$ is a constant. 

In the SA framework, let us have the random vector $\mathbf{X}=(X_{1,}...,X_{d})$ where the variables are mutually independent, and the output $Y=f(\mathbf{X})$ of a deterministic model $f(\cdot)$. 
Thus a functional decomposition
of the variance is available, often referred to as functional ANOVA (Efron and Stein \cite{efrste81}):
\begin{equation}
\mbox{Var}(Y) = \sum_{i=1}^{d}D_{i}(Y)+\sum_{i<j}^{d}D_{ij}(Y)+\cdots+D_{12...d}(Y)
\end{equation}
where $D_{i}(Y)=\mbox{Var}[\mathbb{E}(Y|X_{i})]$, $D_{ij}(Y)=\mbox{Var}[\mathbb{E}(Y|X_{i},X_{j})]-D_{i}(Y)-D_{j}(Y)$
and so on for higher order interactions. The so-called ``Sobol' indices''
or ``variance-based sensitivity indices'' (Sobol \cite{sob93}) are obtained
as follows:

\begin{equation}\label{eqindordre1}
S_{i}=\frac{D_{i}(Y)}{\mbox{Var}(Y)},\quad S_{ij}=\frac{D_{ij}(Y)}{\mbox{Var}(Y)},\quad\cdots
\end{equation}
These indices express the share of variance of $Y$ that is due to
a given input or input combination.

The number of indices growths in an exponential way with the number
$d$ of dimension: there are $2^{d}-1$ indices. For computational
time and interpretation reasons, the practitioner should not estimate
indices of order higher than two. Homma and Saltelli \cite{homsal96}
introduced the so-called ``total indices'' or ``total effects''
that write as follows:
\begin{equation}
S_{T_{i}}=S_{i}+\sum_{i<j}S_{ij}+\sum_{j\neq i,k\neq i,j<k}S_{ijk}+... = \sum_{l\in\#i}S_{l}
\end{equation}
where $\#i$ are all the subsets of $\left\{1,...,d\right\}$ including
$i$. In practice, when $d$ is large, only the main effects and the
total effects are computed, thus giving a good information on the
model sensitivities. 

To estimate Sobol' indices, Monte Carlo sampling based methods have been developed: Sobol \cite{sob93} for first order and interaction indices and  Saltelli \cite{sal02} for fist order and total indices. Unfortunately, to get precise estimates of sensitivity indices, these methods are costly in terms of number of model calls 
 (rate of convergence in $\sqrt{n}$ where $n$ is the sample size). In common practice, the value of $10^{4}$ model calls can be required to estimate the Sobol' index of one input with an uncertainty of $10\%$.
Using quasi-Monte Carlo sequences instead of Monte Carlo samples can sometimes reduce this cost by a factor ten (Saltelli et al. \cite{salrat08}).
The FAST method (Cukier et al. \cite{cuklev78}), based on a multi-dimensional Fourier transform, is  also used to reduce this cost.
Saltelli et al. \cite{saltar99} have extended this technique to compute total Sobol' indices and Tarantola et al. \cite{targat06} have coupled FAST with a  Random Balance
Design. Tissot and Prieur \cite{tispri12} have recently analyzed and improved these methods.
However, FAST remains costly, unstable and biased when the number of inputs increases (larger than $10$) (Tissot and Prieur \cite{tispri12}).

One advantage of using a Monte Carlo based method is that it provides error made on indices estimates via random repetition (Iooss et al. \cite{ioovan06}), asymptotic formulas (Janon et al. \cite{jankle13}) or bootstrap methods (Archer et al. \cite{arcsal97}).
Thus, other Monte Carlo based estimation formulas have been introduced and greatly improved the estimation precision: Mauntz formulas (Sobol et al. \cite{sobtar07}, Saltelli et al. \cite{salann10}) for estimating  small indices, Jansen formula (Jansen \cite{jan99}, Saltelli et al. \cite{salann10}) for estimating total Sobol' indices and Janon-Monod formula (Janon et al. \cite{jankle13}) for estimating large first-order indices.

\begin{figure}[!ht]
\centering
\includegraphics[height=12cm,width=12cm]{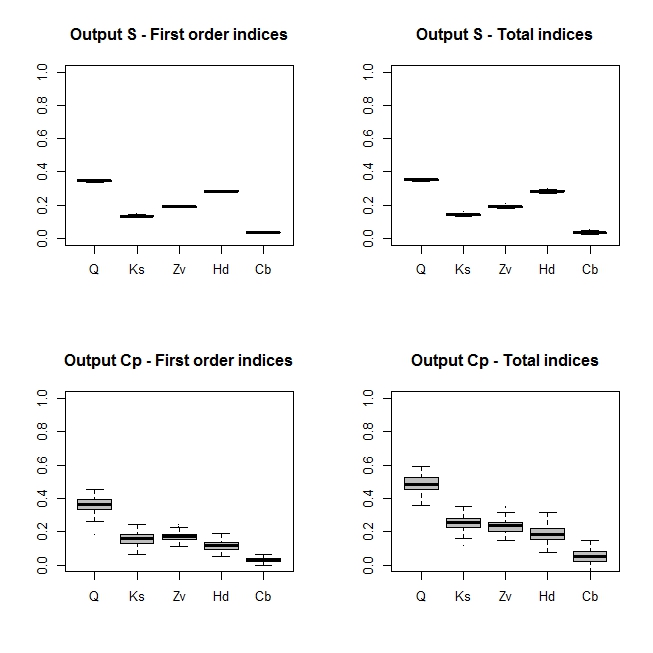}
\caption{Estimation of Sobol' indices on the flood example. Each boxplot corresponds to $100$ independent estimates.}\label{fig:crues_sobol}
\end{figure}

To illustrate the estimation of the Sobol' indices on the flood exercise (Eqs. (\ref{eq:cruesS}) and (\ref{eq:cruesCp})) with $d = 5$ random inputs, we use Saltelli \cite{sal02} formula with a
Monte Carlo sampling. It has a cost $N=n(d+2)$ in terms of model calls where $n$ is the size of an initial Monte Carlo sample.
Here, $n = 10^5$ and we repeat $100$ times the estimation process  to obtain confidence intervals (as boxplots) for each estimated indices. 
Figure \ref{fig:crues_sobol} gives the result of these estimates, which have finally required $N = 7\times 10^7$ model calls.

For the output $S$, the first order indices are almost equal to the total indices, and
results seem very similar to those of SRC$^2$. The model is linear and
the estimation of Sobol' indices is unnecessary in this case. 
For the output $C_p$, we obtain different information than  those provided by SRC$^2$ and SRRC$^2$: the total effect of $Q$ is
about $50\%$ (twice than its SRC$^2$), the effect of $H_d$ is about $20\%$, while $Q$ and $K_s$ have non-negligible interaction effects. Second order Sobol' index between $Q$ and $K_s$ is worth $6\%$.

%%%%%%%%%%%%%
\subsection{Other measures}\label{sec:other}

From an independent and identically distributed sample (as a Monte Carlo one), other techniques can be used for SA.
For example, statistical testing based techniques consist, for each input, to divide the sample into several sub-samples (dividing the considered input into equiprobable stratas).
Statistical tests can then be applied to measure the homogeneity of populations between classes: common means (CMN) based on a Fisher test, common median (CMD) based on a $\chi^2$-test, common variances (CV) based on a Fisher test, common locations (CL) based on the Kruskal-Wallis test, \ldots (Kleijnen and Helton \cite{klehel99}, Helton et al. \cite{heljoh06}).
These methods do not require assumptions about the monotony of the output with respect to the inputs but lacks of some quantitative interpretation.

The indices of \S \ref{sec:sobol} are based on the second-order moment (i.e. the variance) of the ouput distribution.
In some cases, variance poorly represents the variability of the distribution.
Some authors have then introduced the so-called moment independent importance measures, which do not require any computation of the output moments.
Two kinds of indices have been defined:
\begin{itemize}
\item The entropy-based sensitivity indices (Krzykacz-Hausmann \cite{krz01}, Liu et al. \cite{liuche06}, Auder and Iooss \cite{audioo08}),
\item The distribution based sensitivity indices (Borgonovo \cite{bor07}, Borgonovo et al. \cite{borcas11}) which consider a distance or a divergence between the output distribution and the output distribution conditionally to one or several inputs.
\end{itemize}
It has been shown that these indices can provide complementary information than Sobol' indices.
However, some difficulties arise in their estimation procedure.

%%%%%%%%%%%%%%%%%%%%%%%
\section{Deep exploration of sensitivities}

In this section, the discussed methods provide additional sensitivity information than just scalar indices.
Moreover, for industrial computer codes with a high computational cost (from several tens of minutes to days), the estimation of Sobol' indices, even with sophisticated sampling methods, are often unreachable.
This section also summarizes a class of methods for approximating the numerical model to estimate  Sobol' indices at a low computational cost, while providing a deeper view of the input variables effects.

%%%%%%%%%%%%%
\subsection{Graphical and smoothing techniques}

\begin{figure}[!ht]
\begin{center}
    \includegraphics[width=11.8cm,height=4cm]{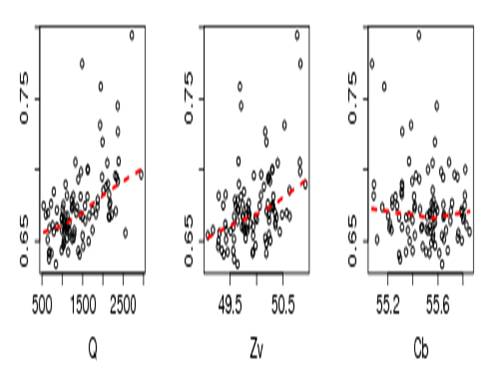} 
    \includegraphics[width=8cm,height=4.5cm]{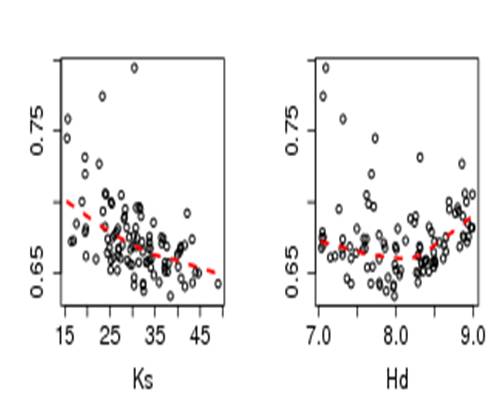} 
\end{center}
\caption{Scatterplots on the flood example with the $5$ inputs $Q$, $Z_v$, $C_b$, $K_s$, $H_d$ and the output $C_p$. Dotted curve is a local polynomial based smoother.}\label{fig:crues_scat}
\end{figure}

Beyond Sobol' indices that only give a scalar value for the effect of an input
$X_i$  on the output $Y$, the influence of $X_i$ on $Y$ along its
domain of variation is also of interest. In the literature, it is often referred to as main effects, but to avoid any
confusion with the indices of the first order, it is preferable to talk about main effects visualization (or
graph). The scatterplots (visualization of point cloud
of any sample simulations $(\mathbf{X}^n,\mathbf{Y}^n)$ with the graphs of $Y$ vs. $X_i$, $i = 1,\ldots,d$)
meets this goal, but in a visual  subjective manner. This is
shown in Figure \ref{fig:crues_scat}, on the flood example and using the $100$-size sample of \S \ref{sec:regr}.

Based on parametric or non-parametric regression methods (Hastie and Tibshirani \cite{hastib90}), the
smoothing techniques aim to estimate the conditional moments of $Y$ at first or higher order. 
SA is often restricted to the determination of the conditional expectation at first and second orders (Santner et al. \cite{sanwil03}), in order to obtain:
\begin{itemize}
\item main effects graphs, between $X_i$ and $\mathbb{E}(Y|X_i)-\mathbb{E}(Y)$ on the whole variation domain of  $X_i$ for $i=1,\ldots,d$;
\\
\item interaction effects graphs, between $(X_i,X_j)$ and $\mathbb{E}(Y|X_iX_j)-\mathbb{E}(Y|X_i)-\mathbb{E}(Y|X_j)-\mathbb{E}(Y)$ on all the variation domain of  $(X_i,X_j)$ for
$i=1,\ldots,d$ and $j=i+1,\ldots,d$.
\end{itemize}

Storlie and Helton \cite{stohel07} conducted a fairly comprehensive review of non-parametric smoothing methods
that can be used for SA: moving averages, kernel methods, local polynomials,
smoothing splines, etc. In Figure \ref{fig:crues_scat}, the local polynomial smoothier is plotted for each cloud of
points, thereby clearly identifying the mean trend of the output versus
each input.

Once these conditional expectations are modeled, it is easy to quantify their variance by sampling, and thus to estimate Sobol' indices (cf. Eq. (\ref{eqindordre1})) of order
one, two, or even higher orders. Da Veiga et al. \cite{davwah08} discuss the theoretical properties of
local polynomial estimators of the conditional expectation and variance, and then deduce
the theoretical properties of the Sobol' indices estimated by local polynomials. 
Storlie and Helton \cite{stohel07} also discuss the efficiency of additive models
 and regression trees to  non-parametrically estimate $\mathbb{E}(Y|X_1,\ldots,X_d)$.
 This finally leads to build an approximate model of $f(\cdot)$, which is called a ``metamodel''.
 This will be detailed in the following section.

\begin{figure}[ht!]
\centering
\includegraphics[width=0.6\textwidth]{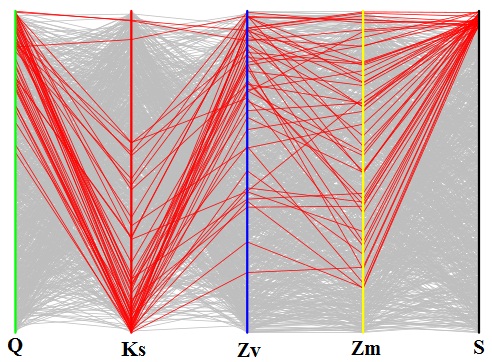}
\caption{Cobweb plot of $10000$ simulations of the flood model.\label{FIGcobWeb}}
\end{figure}

In SA, graphical techniques can also be useful.
For example, all the scatterplots between each input variable and the model output can detect some trends in their functional relation (see Figure \ref{fig:crues_scat}).
However scatterplots do not capture some interaction effects between the inputs.
Cobweb plots (Kurowicka and Cooke \cite{kurcoo06}), also called parallel coordinate plots, can be used to visualize the simulations as a set of trajectories.
In Figure \ref{FIGcobWeb}, the simulations leading to the $5\%$ largest values of the model output $S$ have been highlighted.
This allows to immediately understand that these simulations correspond to large values of the flowrate $Q$ and small values of the Strickler coefficient $K_s$.

%%%%%%%%%%%%%
\subsection{Metamodel-based methods}

%The metamodel is a synonym of response surface, emulator, surrogate model and proxy.
The metamodel concept is frequently used
to simulate the behavior of an experimental system or
a long running computational code based on a certain
number of output values. Under the name of the response surface methodology, it was originally proposed as
a statistical tool, to find the operating conditions of a
process at which some responses were optimized (Box and Draper \cite{boxdra87}). Subsequent
generalizations led to these methods being used to
develop approximating functions of deterministic
computer codes (Downing et al. \cite{dowgar85}, Sacks et al.
\cite{sacwel89}, Kleijnen and Sargent \cite{klesar00}). It consists in generating
a surrogate model that fits the initial data (using for example a least squares procedure), which
has good prediction capabilities and has negligible computational cost.
It is thus efficient for uncertainty
and SA requiring several thousands
of model calculations (Iooss et al. \cite{ioovan06}).

In practice, we focus on three main issues during the construction of a
metamodel:
\begin{itemize}
\item the choice of the metamodel that can be derived from any linear regression model, non-linear parametric or non-parametric (Hastie et al. \cite{hastib02}). 
The most used metamodels include polynomials,
splines, generalized linear models, generalized additive models, kriging, neural networks, SVM, boosting regression trees (Simpson
et al. \cite{simpep01}, Fang et al. \cite{fanli06}). 
Linear and quadratic
functions are commonly considered as a first
iteration. Knowledge on some input interaction types
may be also introduced in polynomials (Jourdan and
Zabalza-Mezghani \cite {jouzab04}, Kleijnen \cite{kle05}). However, these kinds of models are not always
efficient, especially in simulation of complex and non-linear phenomena.
For such models, modern statistical learning
algorithms can show much better ability to build
accurate models with strong predictive capabilities
(Marrel et al. \cite{marioo08});
\item the design of (numerical) experiments. The main qualities required for an experimental design
are the robustness (ability to analyze different models), the effectiveness
(optimization of a criterion), the goodness of points repartition (space filling property) and the low cost for its construction (Santner et al. \cite{sanwil03},
Fang et al. \cite{fanli06}). Several studies have shown the qualities of different types of experimental designs with respect to the predictivity metamodel (for example Simpson et al \cite{simpep01});
\item the validation of the metamodel. In the field of classical experimental design,
proper validation of a response surface is a crucial aspect and is considered with care. 
However, in the field of numerical experiments,
this issue has not been deeply studied. The usual practice is to estimate
global criteria (RMSE, absolute error, ...) on
a test basis, via cross-validation or bootstrap (Kleijnen and Sargent \cite{klesar00},
Fang et al. \cite{fanli06}). When the number of calculations is small and to overcome
problems induced by the cross validation process, Iooss et al. \cite{ioobou10} have recently studied how to minimize the size of a test sample, while obtaining a
good estimate of the metamodel predictivity.
\end{itemize}

Some metamodel allows to directly obtain the sensitivity indices. For example,
Sudret \cite{sud08} has shown that Sobol' indices  are a by-product result
of the polynomial chaos decomposition. 
The formulation of the kriging metamodel provides also analytical formula for the Sobol' indices, associated with interval confidence coming from the kriging error (Oakley and O'Hagan \cite{oakoha04}, Marrel et al. \cite{marioo09}, Le Gratiet et al. \cite{legcan13}).
A simplest idea, widely used in practice, is to apply an intensive sampling technique (see \S \ref{sec:sobol}) directly on the metamodel
to estimate Sobol' indices (Santner et al. \cite{sanwil03}, Iooss et al. \cite{ioovan06}). The variance proportion  not
explained by the metamodel (calculated by $1-Q_2$, cf. Eq. (\ref{eqQ2})) gives us what is missing in the SA (Sobol \cite{sob03}). Storlie et al. \cite{stoswi09} propose a
bootstrap method for estimating the impact of the metamodel error.

As in the previous section, we can be interested by visualizing main effects (Schonlau and Welch \cite{schwel06}). These can be directly given by the metamodel
(this is the case with the polynomial chaos methods, kriging, additive models), or computed by simulating the
conditional expectation $\mathbf{E}(Y|X_i)$.

To illustrate our purpose, we use the flood example (Eqs. (\ref{eq:cruesS}) and (\ref{eq:cruesCp})). A kriging metamodel  is built on a $100$-size Monte Carlo sample with inputs $Q$, $K_s$, $Z_v$, $H_d$, $C_b$ and on the output $C_p$. The metamodel consists in a deterministic term (simple linear model), and a corrective term modeled by a Gaussian stationary stochastic process, with a generalized exponential covariance (see Santner et al \cite{sanwil03} for more details). The technique for estimating the metamodel hyperparameters is
described in Roustant et al. \cite{rougin12}. The predictivity coefficient estimated by leave-one-out is
$Q^2 = 99\%$ compared with $Q^2 = 75\%$ obtained with a simple linear model. The kriging metamodel
 is then used to estimate Sobol' indices  in the same manner as in
\S \ref{sec:sobol}: Saltelli's estimation formula, Monte Carlo sampling, $n = 10^5$, $r = 100$ repetitions.
This requires $N = 7\times 10^7$ metamodel predictions. In Table \ref{tab:sobol}, we compare
Sobol' indices (averaged over $100$ repetitions) obtained with the metamodel to those obtained
with the ``real'' flood model (Eqs. (\ref{eq:cruesS}) and (\ref{eq:cruesCp})). Errors between these two estimates are relatively
low: with only $100$ simulations with the true model, we were able to obtain precise estimates
 (errors $<15\%$) of first order and total Sobol' indices.

\begin{table}[!ht]
  \centering
\caption{Sobol' indices estimated by Monte Carlo sampling (cost of $N=7\times 10^7$ evaluations) using the flood model and a metamodel fitted on $N'=100$ calls of the flood model.}\label{tab:sobol}
\begin{tabular}{lccccc}
Indices (in $\%$) & $Q$ & $K_s$ & $Z_v$ & $H_d$ & $C_b$ \\
\hline\hline
$S_i$ model & 35.5 & 15.9 & 18.3 & 12.5 & 3.8 \\
\hline
$S_i$ metamodel & 38.9 & 16.8 & 18.8 & 13.9 & 3.7 \\
\hline\hline
$S_{T_i}$ model & 48.2 & 25.3 & 22.9 & 18.1 & 3.8 \\
\hline
$S_{T_i}$ metamodl & 45.5 & 21.0 & 21.3 & 16.8 & 4.3
\end{tabular}
\end{table}

%%%%%%%%%%%%%%%%%%%%%%%
\section{Synthesis and conclusion}

Although all SA techniques have not been listed, this review has illustrated the great
 variety of  available methods, positioning in terms of
assumptions and kind of results. 
Moreover, some recent improvements have not been explained, for example for the Morris method (Pujol \cite{puj09}). 

\begin{figure}[!ht]
\centering
\includegraphics[width=\textwidth]{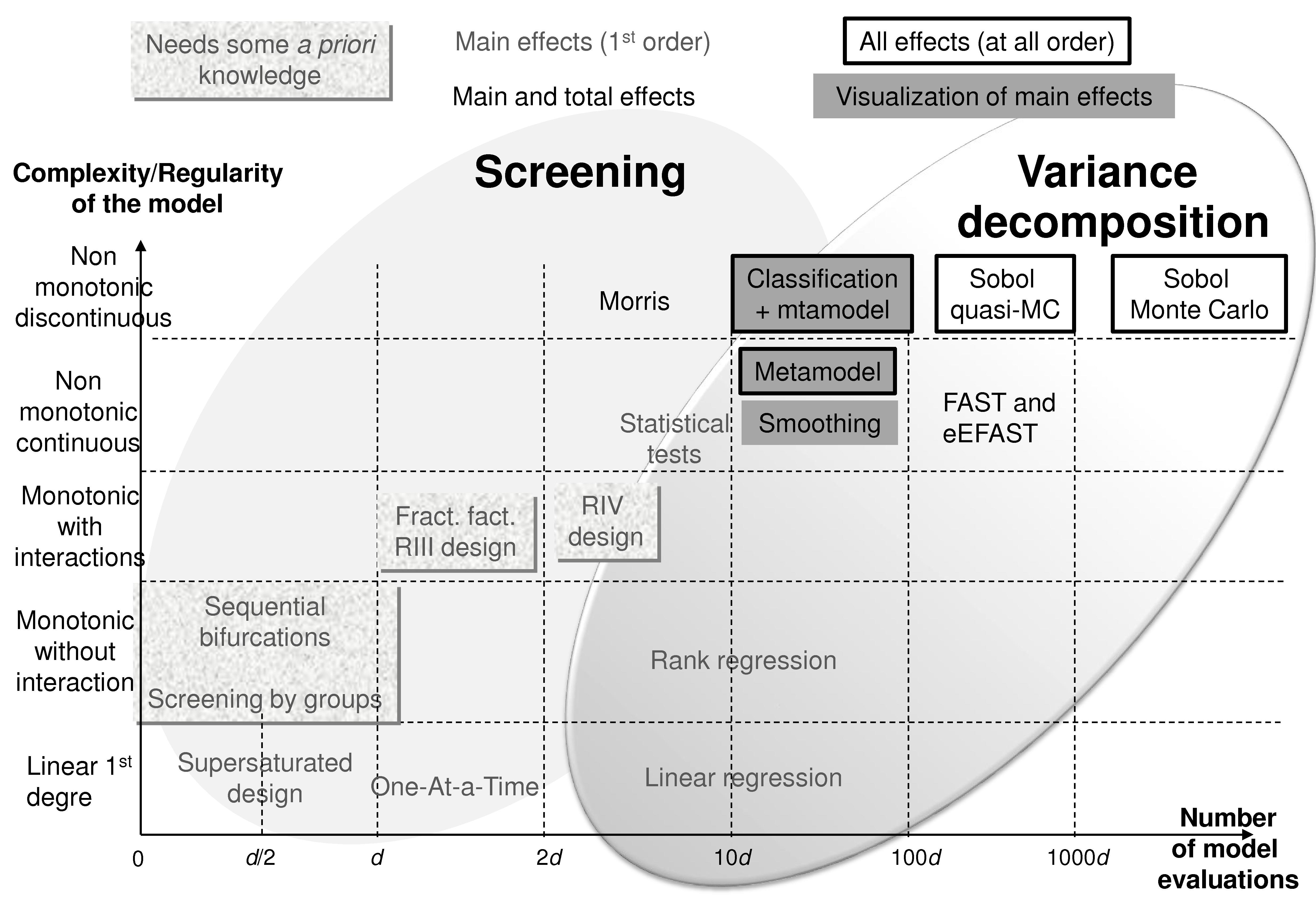}
\caption{SA methods graphical synthesis.}
\label{fig:syntheseGSA2}
\end{figure}

A synthesis is provided in Figure \ref{fig:syntheseGSA2} which has several levels of reading:
\begin{itemize}
\item distinction between screening methods (identification of non-influential variables among a large number) and more precise variance-based quantitative methods,
\item positioning methods based on their cost in terms of model calls number (which linearly depends in the number of inputs for most of the methods),
\item positioning methods based on their assumptions about the model complexity and regularity,
\item distinction between the type of information provided by each method,
\item identification of methods which require some a priori knowledge about the model behaviour.
\end{itemize}

\begin{figure}[!ht]
        \centering
\includegraphics[width=\textwidth]{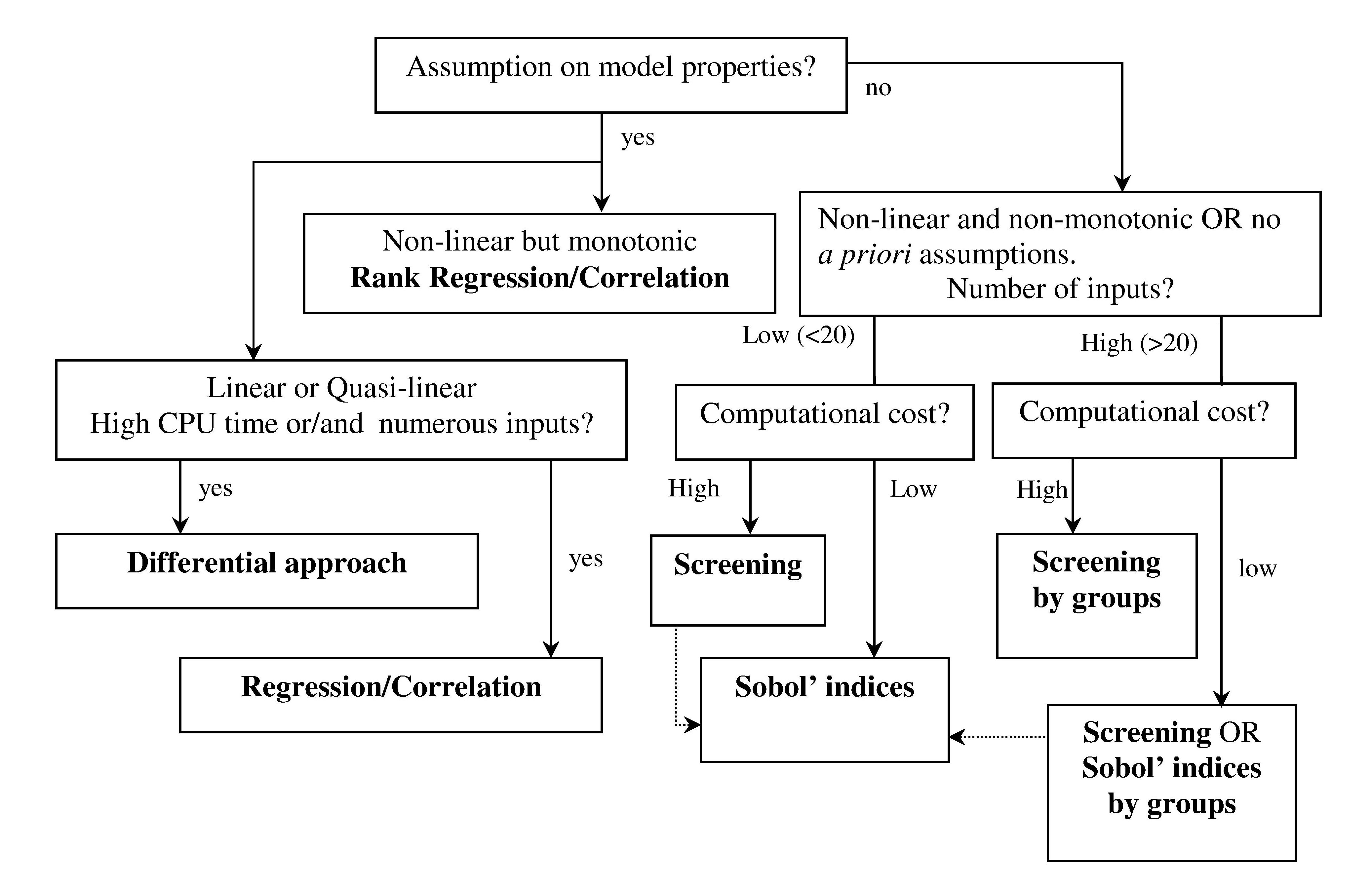}
\caption{Decision diagram for the choixe of a SA method (from de Rocquigny et al. \cite{derdev08}).}\label{fig:decisiondiagram}
\end{figure}

Based on the characteristics of the different methods, some authors (de Rocquigny et al. \cite{derdev08}, Pappenberger et al. \cite{paprat10}) have proposed decision trees to help the practitioner to
choose the most appropriate method for its problem and its model. Figure \ref{fig:decisiondiagram} reproduces the
flowchart of de Rocquigny et al. \cite{derdev08}. Although useful to fix some ideas, such diagrams are rather simple and should be used with caution.

Several issues about SA remain open.
For instance, recent theoretical results have been obtained on the asymptotical properties and efficiency of Sobol' indices estimators (Janon et al. \cite{jankle13}), but estimating total Sobol' indices at low cost is a problem of primary importance in applications (see Saltelli
et al. \cite{salann10} for a recent review on the subject). 
SA for dependent inputs
has also been discussed by several authors (Saltelli and Tarantola \cite{saltar02}, Jacques
et al. \cite{jaclav06}, Xu and Gertner \cite{xuger07}, Da Veiga et al. \cite{davwah08},Gauchi et al. \cite{gauleh10}, Li et al. \cite{lirab10},
Chastaing et al. \cite{chagam12}), but this issue remains misunderstood. 

This chapter has been focused on SA relative to the overall variability of model output. In practice, one can be interested by other quantities of interest, such as the output entropy (cf. \S \ref{sec:other}), the probability that the output
exceeds a threshold (Saltelli et al. \cite{salcha00}, Frey and Patil \cite{frepat02}, Lema\^{\i}tre et al. \cite{lemser13}) or a quantile estimation (Cannamela et al. \cite{cangar08}). 
This is an active area of research.

In many applications, the model output is not a single scalar but a vector or a function (temporal, spatial, spatio-temporal, \ldots).
Campbell et al. \cite{cammck06}, Lamboni et al. \cite{lammon10}, Marrel et al. \cite{marioo11} and Gamboa et al \cite{gamjan13} have
produced first SA results on such problems. The case of
functional inputs also receives a growing interest (Iooss and Ribatet
\cite{ioorib08}, Lilburne and Tarantola \cite{liltar09}, Saint-Geours et al. \cite{sailav11}), but its
treatment in a functional statistical framework remains to be done.

In some situations, the computer code is not a  deterministic simulator but a stochastic one. This means that two model calls with the same set of input variables leads to
different output values. 
Typical stochastic computer codes are queuing models, agent-based models,  models involving partial differential equations applied to heterogeneous or Monte-Carlo based numerical models.
 For this type of codes,  Marrel et al. \cite{marioo12} have proposed a first solution for dealing with Sobol' indices.

Finally, quantitative SA methods are limited to low-dimensional models, with no more than a few tens of input variables.
On the other hand, deterministic methods, such as adjoint-based ones (Cacuci \cite{cac03}), are well suited when the model includes a large number of input variables.
A natural idea is to use the advantages of both methods.
Recently introduced, Derivative-Based Sensitivity Measures (DGSM), consists in computing the integral of the square model derivatives for each input (Sobol and Kucherenko \cite{sobkuc09}).
An inequality relation has been proved between total Sobol' indices and DGSM which allow to propose some interpretative results (Lamboni et al. \cite{lamioo13}, Roustant et al. \cite{roufru13}).
It opens the way to perform global SA in high dimensional context.

%%%%%%%%%%%%%%%%%%%%%%%
\section{Acknowledgments}
Part of this work has been backed by French National Research
Agency (ANR) through COSINUS program (project COSTA BRAVA
no. ANR-09-COSI-015).
We thank Anne-Laure Popelin for providing the cobweb plot.

%%%%%%%%%%%%%%%%%%%%%%%%%%%

%\bibliographystyle{plain}
%\bibliography{BIBLIOGRAPHIE/books,BIBLIOGRAPHIE/articles,BIBLIOGRAPHIE/inbooks,BIBLIOGRAPHIE/nt,BIBLIOGRAPHIE/thesis,BIBLIOGRAPHIE/proceedings,BIBLIOGRAPHIE/rapports,BIBLIOGRAPHIE/soft,BIBLIOGRAPHIE/presentation}

\end{document}